\documentclass{amsart}

\usepackage{amsthm}
\usepackage{color, xcolor}

\usepackage{graphicx,epsfig}
\usepackage{epstopdf}
\usepackage{amsmath, amssymb,float, latexsym, euscript,picins}
\usepackage{caption, wrapfig, multirow, tabularx, mathrsfs,verbatim}

\usepackage{url}
\usepackage[all]{xy}
\usepackage{psfrag}
\usepackage{placeins}
\usepackage{float}

\definecolor{ivory}{rgb}{1.0, 1.0, 0.94}
\definecolor{cadmiumgreen}{rgb}{0.0, 0.42, 0.24}
\definecolor{antiquewhite}{rgb}{0.98, 0.92, 0.84}

\def\RR{\mathbb R}
\def\HH{\mathbb H}

\def\SL{\operatorname{SL}}

\def\Ncal{\mathcal N}

\def\Teich{\mathcal T}

\def\MF{\mathcal M\mathcal F}

\def\bdy{\partial}

\def\Tcal{\mathcal T}

\def\MF{\mathcal M\mathcal F}
\def\Ccal{\Delta}

\def\bdy{\partial}

\def\Hcal{\mathcal H}

\def\Acal{\mathcal A}
\def\Shorts{{\bf SH}}
\def\Annuli{{\bf AN}}
\def\Dcal{\mathcal D}

\def\collardepth{{\bf \kappa}}

\def\hyp{\operatorname{{\bf hyp}}}
\def\mf{\operatorname{{\bf mf}}}

\def\Scal{\mathcal S}

\def\ZZ{\mathbb Z}

\def\Lcal{\mathcal L}

\def\stdan{\mathbb A}

\def\tr{\operatorname{tr}}

\def\Qcal{\mathcal Q}

\newcommand{\R}{\mathbb R}

\theoremstyle{definition}

\newtheorem{theorem}{Theorem}[section]
\newtheorem{lemma}[theorem]{Lemma}

\newtheorem{definition}[theorem]{Definition}
\newtheorem{proposition}[theorem]{Proposition}

\newtheorem*{mainthmA}{Main Theorem A}
\newtheorem*{mainthmB}{Main Theorem B}
\newtheorem*{mainthmC}{Main Theorem C}

\begin{document}

\pagecolor{ivory}

\title[Collar parameters for Teichm\"uller Space \& Measured Foliations on a Surface]{Collar parameters for Teichm\"uller Space \& Measured Foliations on a Surface Research Announcement}
\author{Daryl Cooper and Catherine Pfaff}

\begin{abstract} We defined a new set of coordinates with respect to which the Thurston compactification of Teichm\"uller space is the radial compactification of Euclidean space. \end{abstract}
 
\maketitle

The seminal work \cite{WPTsurfaces} of Thurston uses lengths of simple closed curves on a surface to define a compactification of its Teichm\"uller space. Let $\Scal$ denote the set of isotopy classes of essential simple closed curves on a closed orientable surface $\Sigma$ of genus $g\ge 2$. Throughout this paper a {\em measured foliation} is a {\em transversally measured singular foliation on a surface}. A hyperbolic metric and a measured foliation on $\Sigma$ each assign a length to members of $\Scal$. Both then determine a projectivized length function on $\Scal$, leading to Thurston's famous \cite{WPTsurfaces} compactification of Teichm\"uller space $\Teich(\Sigma)$.

\parpic[r]{\includegraphics[width=2.6in]{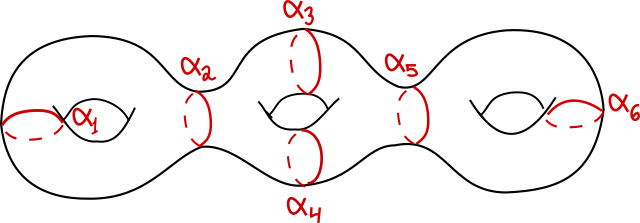}}
As depicted to the right, let $\{\alpha_i\}$ denote a set of pairwise disjoint simple closed curves on $\Sigma$ whose complement is a disjoint union of 3-holed spheres (pants). Fenchel-Nielsen coordinates on Teichm\"uller space assign a length and twist parameter to each $\alpha_i$. Dehn-Thurston coordinates use similar data to parameterize measured foliations. In both cases the length is positive, or possibly zero in the case of  measured foliations, but the twist is an arbitrary real number. The {\em collar parameter (CP) coordinates} we define here assign a point in $\RR^2$ to each $\alpha_i$. They are a variant of the Fenchel-Nielsen and Dehn-Thurston coordinates, and encode both the length and the twist parameter. With respect to CP coordinates the Thurston compactification is the radial compactification of Teichm\"uller space. 

We now describe CP coordinates. In what follows the term {\em structure} means either a hyperbolic metric or a measured foliation on $\Sigma$. 
Each $\alpha_i$ is contained in an annulus $A_i$ that in some sense is maximal in the structure. In a hyperbolic structure, $A_i$ is provided by the
{\em Collar Lemma} \cite{Buser}. By Proposition \ref{standardfolation} every measured foliation is equivalent to one such that either $\bdy A_i$ is transverse to the foliation, or else $A_i$ is a union of smooth closed leaves such that the union of the closed leaves isotopic into $A_i$, that are in the complement of $A_i$, has zero transverse measure.
 Such an annulus in either structure is called a {\em standard collar}.

Given a structure on $\Sigma$ and pair of pants in a pants decomposition, there is a structure-preserving reflection of the pants whose fixed point set consists of three arcs, one arc connecting each pair of boundary components. The arcs are called {\em seams}. Each annulus $A_i$ has two basepoints on each boundary component given by the intersection with the seams of the pants decomposition.
\parpic[l]{\includegraphics[width=1.85in]{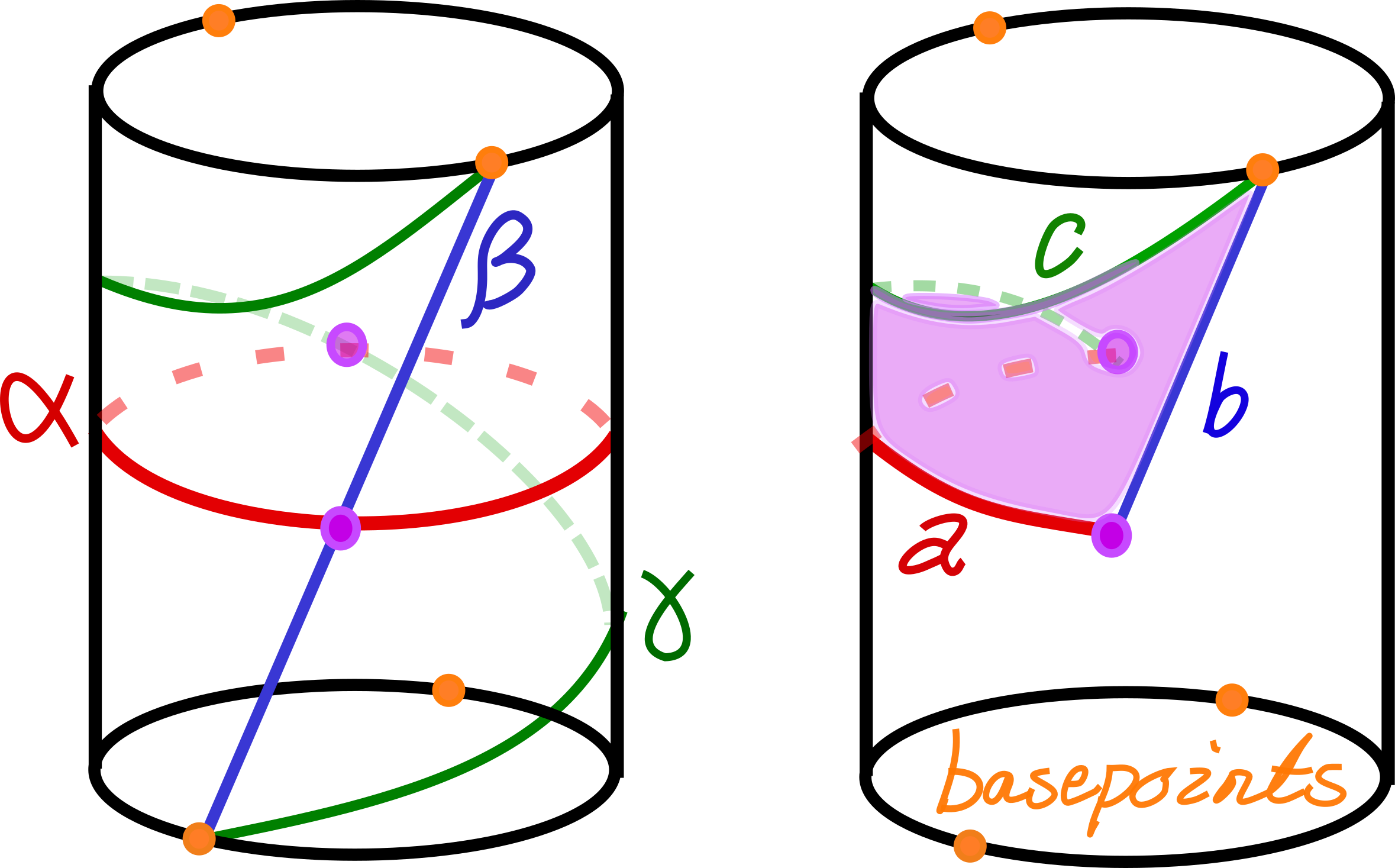}}
\noindent  CP coordinates parameterize structures on the annulus up to isotopy fixing the basepoints. 
Triangulate the annulus using these basepoints, and arcs $\beta$ and $\gamma$ connecting them, plus three meridian circles around the annulus. The structure on the annulus is determined by the lengths of the sides  of a triangle in this triangulation ($a$,$b$,$c$ in the  image). There is an equation relating these lengths, giving a parameter space $\RR^2$ with coordinates that are certain linear combinations of edge lengths. This parameter space is called the space of {\em collar parameters}, and is described in \S \ref{ss:CPcoords}. It is a pleasant fact that for both structures the equation is symmetric in the three edge lengths.

\parpic[r]{\includegraphics[width=1.8in]{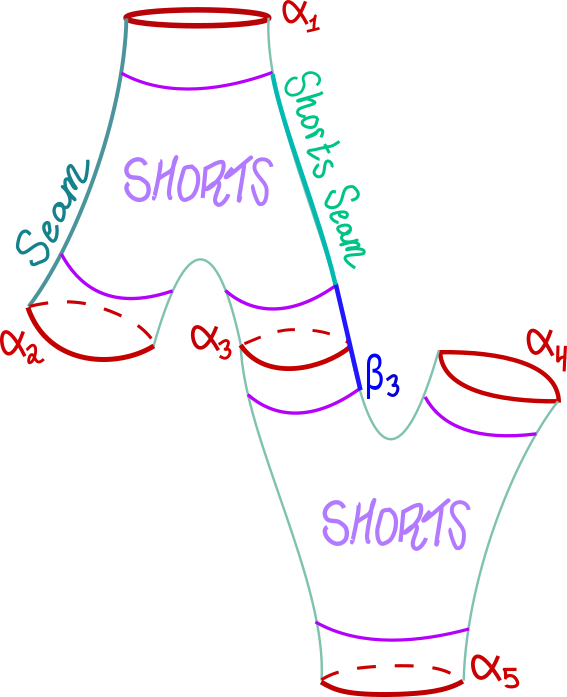}} 
In the hyperbolic case, the curves $\{\alpha_i\}$ are geodesics and separate $\Sigma$ into hyperbolic pants.
Deleting the interiors of the annuli results in subsurfaces called {\em shorts}. The boundary components of the shorts are hypercycles (curves that are equidistant from a geodesic). Every point in the shorts is within distance $\operatorname{cosh}^{-1}(3)$ of the boundary. 

The intersection of a seam of the pants with the shorts  is the unique geodesic arc in the shorts connecting that pair of boundary components and meeting them orthogonally.  These arcs are called the {\em seams} of the shorts. The shorts are determined up to isometry by the data for the annuli. Gluing the shorts to the annuli so that the basepoints on the annuli are endpoints of seams parameterizes Teichm\"uller space by a product of parameter spaces for the annuli. 
 
A similar procedure works for measured foliations. Each measured foliation on the annulus is {\em linear} (see \S \ref{ss:Foliations}). The measured foliations on the pairs of shorts are determined solely  by the measure of the boundary components. Again there are seams: the fixed points of a reflection that preserves the measure. These foliations yield a global parametrization of the space of measured foliations on $\Sigma$, as a product of the parameter spaces for the annuli.

This gives the parameterizations $\Theta_{\Tcal}: \RR^{6g-6}\to\Teich(\Sigma)$ of Teichm\"uller space, and $\Theta_{\MF}:\RR^{6g-6}\to\MF(\Sigma)$ of the space of measured foliations, using CP coordinates. 
Since both sets of coordinates are determined by the lengths of the {\em same sides of the same triangles} the homeomorphism $\Theta_{\Tcal}\circ \Theta_{\MF}^{-1}$, sending a  measured foliation to a hyperbolic metric, is a {\em good approximation} for large foliations. This works so well because most of the length of a geodesic (after a small perturbation), and all the measure, is concentrated in the collars. Using these coordinates on Teichm\"uller space, it follows that the Thurston compactification is just the radial compactification of Euclidean space. The same result does not hold using Fenchel-Nielsen coordinates, as is apparent by considering sequences where the lengths of some $\alpha_i$ go to zero. There are explicit formulae (Proposition \ref{FNDT}) for collar parameters in terms of the Fenchel-Nielsen coordinates or the Dehn-Thurston coordinates.

\parpic[l]{\includegraphics[width=1.3in]{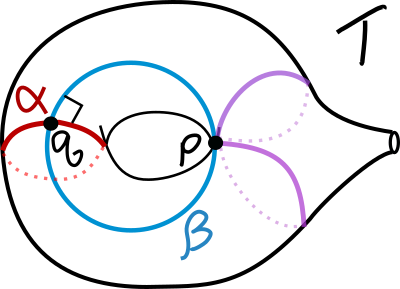}} 
A hyperbolic torus with a geodesic boundary component can be obtained by identifying two 
boundary components of a hyperbolic pair of pants that have the same length ($\alpha$ in the image). The limit, as the length of
the torus' boundary component goes to zero, is a complete hyperbolic once-punctured torus $T$ with finite area. Given a standard collar $\stdan$ containing a closed geodesic $\alpha$, there is such a $T$ that contains an isometric image of the interior of $\stdan$ and exactly  two points, one on each (purple) boundary component of $\stdan$, are identified to a single point $p$ in $T$. Twisting along $\alpha$ produces a one-parameter family of such tori. A geodesic arc crossing $\stdan$ and with these endpoints gives a  geodesic loop $\beta$ in $T$. There is another loop $\gamma$ such that $\alpha,\beta,\gamma$ are three closed geodesics in $T$ that pairwise-intersect precisely once. With suitable orientations, we have $\alpha\cdot\beta\cdot\gamma=1\in\pi_1T$. These closed geodesics contain the three edges in a triangle of the triangulation of $\stdan$, and hence define the collar parameters. More precisely, the length of a triangle edge is half the length of the geodesic loop that contains it. The commutator $[\alpha,\beta]$ is parabolic. The formula for the trace of this parabolic, expressed in terms of the lengths of $\alpha,\beta,\gamma$, gives the equation relating the edge lengths in a standard collar. {\em It follows that one may regard the collar parameters of a standard collar as a point in the Teichm\"uller space of finite area complete hyperbolic metrics on $T$, and thus can view $\Teich(\Sigma)$ as a product of $(3g-3)$ copies of the Teichm\"uller space of the punctured torus.}

Given $v\in\RR^{6g-6}$ with $|v|\le 1$ one can write down an explicit quadratic differential on $\Sigma$ that varies continuously with $v$.
When $|v|<1$ this is a rescaling of a hyperbolic metric with collar parameters $v/(1-|v|)$. For $|v|=1$ it is a measured foliation for collar parameter $v$. This realizes the Thurston compactification of Teichm\"uller as a subspace of $\Qcal(\Sigma)$, the space of quadratic differentials on $\Sigma$.

Theorem A bounds the difference between the length of an isotopy class of a loop in a hyperbolic metric and the corresponding measured foliation. The bound is in terms of the minimum word length in the conjugacy class for the loop. 
This implies Theorem B concerning the compactification. These theorems follow from a stronger result that compares pointwise the hyperbolic metric and corresponding measured foliation
after isotoping these structures into a nice position. 
This culminates in Theorem C, which lifts both Teichm\"uller space and the space of measured foliations to spaces of quadratic differentials where the Thurston compactification arises from (rescaled) quadratic differentials, rather than isotopy classes of structures.

\vspace{5mm}

%%%%%%%%%%%%%%%%%%%%%%%%%%%%%%%%%%%%%%%%%%%%%%%%%%%%%%%%%%%%%%%
 
\section{Collar parameters}\label{s:CollarParameters}

%%%%%%%%%%%%%%%%%%%%%%%%%%%%%%%%%%%%%%%%%%%%%%%%%%%%%%%%%%%%%%%

In what follows, $\Sigma$ is  a closed orientable surface of genus $g\ge 2$ and a {\em structure} is either a hyperbolic metric or a measured foliation on $\Sigma$. 

\begin{figure}[H]
\includegraphics[height=1.3in]{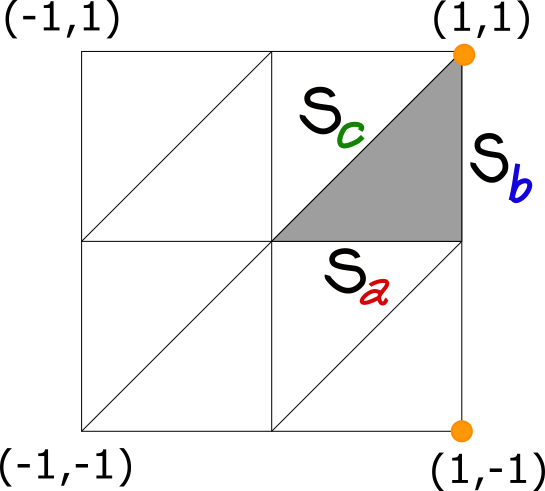}
        \hspace{15mm}
\includegraphics[height=1.27in]{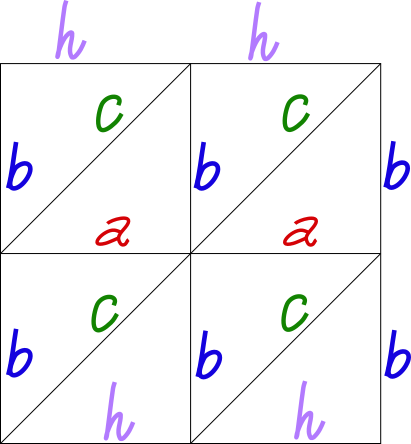}
\vspace*{-1mm}
     \caption{(a) The left-hand image depicts a triangulation of a fundamental domain for the standard annulus $\stdan = S^1\times[-1,1]$ with the reference triangle shaded in. The basepoints are in orange. (b) The right-hand image shows the lengths of the edges determined by a hyperbolic metric or measured foliation on $\stdan$.}
    {\label{fig:Triangulation}}
\end{figure}

To provide a common frame of reference for the structures, we fix a triangulation of a {\em standard annulus}. Define the circle $S^1:=\RR/2\ZZ$, and the {\em standard annulus} as $\stdan := S^1\times[-1,1]$. The universal cover of the standard annulus is   $\widetilde{\stdan}=\RR\times[-1,1]\subset\RR^2$. Let $p :\RR\times[-1,1]\rightarrow \stdan$ be the  covering map. The subset $[-1,1]^2$ of $\widetilde\stdan$ is a fundamental domain. Triangulate  $[-1,1]^2$ with eight Euclidean triangles as shown in Figure \ref{fig:Triangulation}a. Their images under $p$ give a triangulation of $\stdan$. This triangulation contains a {\em reference triangle} with sides that are $s_a=p([0,1]\times 0)$ and $s_b=p(1\times[0,1])$ and $s_c=p(\{(t,t):0\le t\le1\})$. The basepoints on $\bdy \stdan$ are $p(1,\pm1)$.

\vspace{-1mm}

In the following we explain how a triple $(a,b,c)$ in $\Hcal$ (see Equation \ref{Heqtn}) or $\Delta$ (see Equation \ref{e:TriangleEquality}) determines, respectively, a hyperbolic metric or a measured foliation on the reference triangle where the lengths of $s_a$, $s_b$, and $s_c$ are respectively $a$, $b$, and $c$. Using the triangulation of Figure \ref{fig:Triangulation}, we then define, respectively, a hyperbolic metric or measured foliation on the standard annulus.

\vspace{5mm}

%%%%%%%%%%%%%%%%%%%%%%%%%%%%%%%%%%%%%%%%%%%%%%%%%%%%%%%%%%%%%%%

\subsection{Shorts decomposition}\label{ss:ShortsDecomp}

%%%%%%%%%%%%%%%%%%%%%%%%%%%%%%%%%%%%%%%%%%%%%%%%%%%%%%%%%%%%%%%

Suppose
$$\Acal=\{\alpha_i:1\le i\le 3g-3\}$$  
\parpic[r]{\includegraphics[width=1.75in]{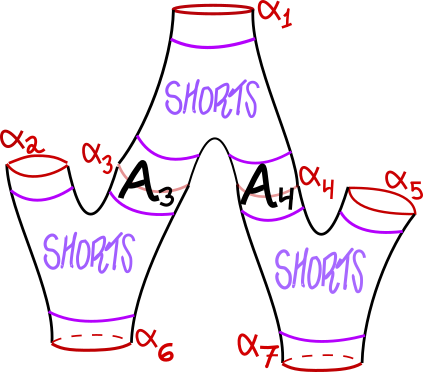}} 
\noindent is a set of disjoint simple closed curves in $\Sigma$ such that the closure of each component of $\Sigma\setminus\Acal$ is a pair of pants. Then there are $2g-2$ complementary components. 

Let ${\Annuli}=\{A_i:1\le i\le 3g-3\}$ be a set of pairwise disjoint compact annuli in $\Sigma$ such that $A_i$ is a neighborhood  of $\alpha_i $.  Let $\Shorts$ be the closure
of $\Sigma\setminus\Annuli$. Then $\Shorts=\{S_j: 1\le j\le 2g-2\}$ where each $S_j$ is called a {\em pair of shorts}, and is a pair of pants with annuli neighborhoods of the boundary components removed.

Then $\Dcal=(\Shorts,\Annuli)$ is called a {\em shorts decomposition} of $\Sigma$.

%%%%%%%%%%%%%%%%%%%%%%%%%%%%%%%%%%%%%%%%%%%%%%%%%%%%%%%%%%%%%%%

\subsection{Hyperbolic collars}\label{ss:HyperbolicCollars}

%%%%%%%%%%%%%%%%%%%%%%%%%%%%%%%%%%%%%%%%%%%%%%%%%%%%%%%%%%%%%%%

\begin{definition}[collars] 
A {\em hyperbolic collar} is an annulus $A$ endowed with a hyperbolic metric so  that $A$ contains a simple closed
geodesic $\alpha$ and each point of $\bdy A$ is a fixed distance $d$ from $\alpha$. Thus $\bdy A$ consists
of two {\em hypercycles}. The 
number $d$ is called the {\em depth} of the collar.
A hyperbolic collar whose core geodesic has length $2a$ is a {\em standard collar} if it has depth $\collardepth(a)=\sinh^{-1}(1/\sinh(a))$.
\end{definition}

The {\em Collar Lemma},  see \cite{Buser}, says that disjoint closed geodesics in a hyperbolic surface are contained in disjoint standard collars.

\vspace{-7mm}

  \begin{center}
\begin{figure}[H]
\includegraphics[height=.8in]{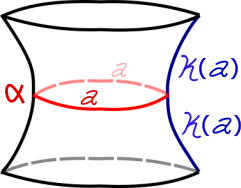}
        \hspace{10mm}
\includegraphics[height=.8in]{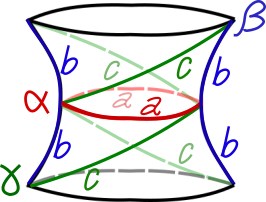}
\vspace*{1mm}
 \caption{(a) A standard collar and (b) its triangulation.} \label{f:standardcollar}
\end{figure}
\end{center}

\vspace{-12mm}

\noindent A hyperbolic collar can be triangulated, as described above,  in such a way so that the sides of the reference triangle are geodesics. 
This is shown in Figure \ref{f:standardcollar}b, however in general $\beta$ is {\em not orthogonal} to $\alpha$. A hyperbolic collar is standard if and only if the side lengths of the reference triangle in the triangulation satisfy the {\em collar equation}:

\begin{lemma}[collar equation]\label{l:boundaryequation} 
If $S$ is a standard collar and the length of the core geodesic is $2a$ then  the edge lengths $(a,b,c)$ of the reference triangle satisfy
$$\cosh^2a+\cosh^2b+\cosh^2c = 2\cosh a\cosh b\cosh c.$$
\end{lemma}

\vspace{3mm}

We say a point in $\RR^3$ satisfies the {\em collar equation} if it is in the set
\begin{equation}\label{Heqtn}
\Hcal=\{(a,b,c)\ :\ \cosh^2a+\cosh^2b+\cosh^2c = 2\cosh a\cosh b\cosh c\quad a,b,c> 0\}.
\end{equation}

\begin{center}
\begin{figure}[H]
\includegraphics[scale=.5]{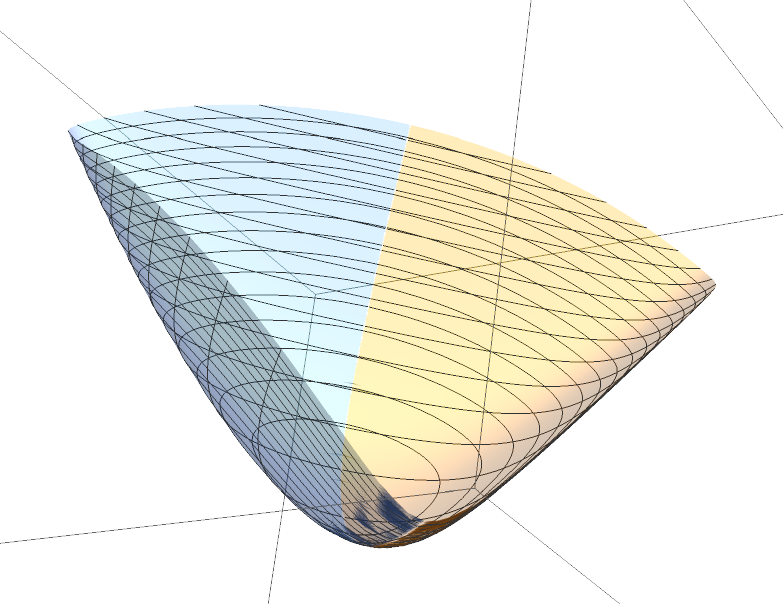}
\hspace{.3in}\includegraphics[scale=.5]{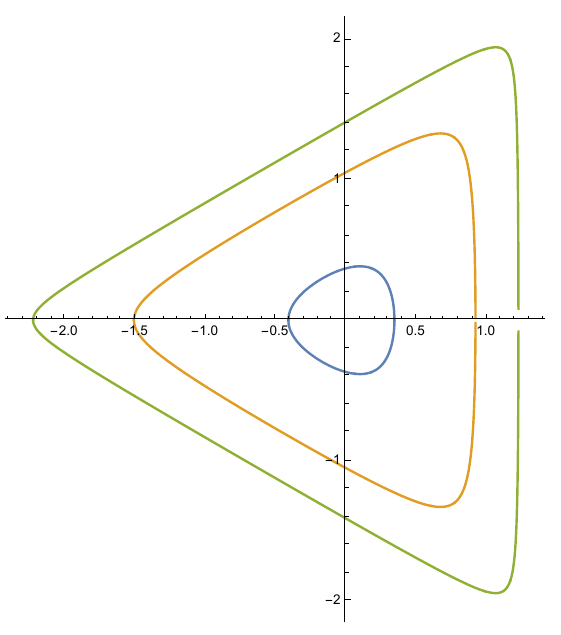}
\caption{(a) On the left is $\Hcal$ and (b) on the right is some cross sections of $\Hcal$} \label{MagicPic}
\end{figure}
\end{center}

The set $\Hcal$ is the left-hand image in Figure \ref{MagicPic}. It sits inside the cone from the origin on the triangle in the plane $x+y+z=2$ with vertices the points $(1,1,0)$ and $(1,0,1)$ and $(0,1,1)$. It is asymptotic to the sides of this cone. The intersections of $\Hcal$ with the planes $x+y+z=C$ are convex curves becoming larger and more nearly triangular as $C$ increases. This is depicted on the right for $C=3,4,5$.
 
\vspace{-6mm}
 
Given $(a,b,c)\in\Hcal$ there is a unique hyperbolic metric $\hyp(a,b,c)$ on $\stdan=S^1\times[-1,1]$ satisfying all of the following properties. 
\begin{itemize}
\item[1.]  The metric is isometric to a standard collar with core curve of length $2a$.
\item[2.]  The metric is preserved by rotation of the $S^1$ factor.
\item[3.] The curve $\beta:[-1,1]\rightarrow\stdan$ given by $\beta(t)=p(-1,t)$ is a geodesic of speed $b$ and length $2b$.
\item[4.] There is a geodesic $\gamma:[-1,1]\rightarrow\stdan$  of length $2c$ and homotopic rel endpoints  to $t\mapsto p(t,t)$.\end{itemize}
 
Rotational invariance implies that each $S^1\times y$ is a hypercycle and has constant speed. Half the length of each boundary component of $\stdan$ is $h=a \coth(a)$.
The sides $s_a$ and $s_b$ of the reference triangle are hyperbolic geodesics, but the third side $s_c$ is not.
The {\em geodesic reference triangle} is the triangle in $\stdan$ with sides $s_a,s_b$ and $p\circ\gamma[0,1]$. It has geodesic
sides and is isotopic to the reference triangle without moving the vertices. It is clear that a standard triangulation of a standard collar  can be isotoped without moving the basepoints to be $\hyp(a,b,c)$. We call such a metric {\em standard}. The next result implies that a standard metric on an annulus differs pointwise from
a linear measured foliation by at most twice the Euclidean metric.

\begin{lemma}\label{collarmetric} The metric $\hyp (a,b,c)$ on $\stdan$ pulls back using the covering space projection to the metric on $\widetilde{\stdan}=\RR\times[-1,1]$ given
 by $$ds^2=\left(a^2+\left(\frac{a}{\sinh a}\right)^2\left(\frac{\sinh(by)}{\sinh b}\right)^2\right) dx^2 \pm 
2ab\sqrt{1-\left(\frac{1}{\sinh a\sinh b}\right)^2} dx dy + b^2 dy^2$$
The sign is $+1$ if $\cosh c\ge\cosh a \cosh b$.
Hence $$|ds^2-(a.dx\pm b.dy)^2| \le 2(dx^2+dy^2)$$ 
\end{lemma}

\vspace{3mm}

%%%%%%%%%%%%%%%%%%%%%%%%%%%%%%%%%%%%%%%%%%%%%%%%%%%%%%%%%%%%%%%%%

\subsection{Triangle lengths and collar parameters}\label{ss:CPcoords}

%%%%%%%%%%%%%%%%%%%%%%%%%%%%%%%%%%%%%%%%%%%%%%%%%%%%%%%%%%%%%%%%%

Let $\pi:\RR^3\rightarrow\RR^2$ be the linear map defined by
$$\pi(a,b,c)=(4a-2b-2c,2b-2c).$$
This is the composition of orthogonal projection of $\RR^3$ onto the subspace given by $a+b+c=0$, followed by an isomorphism to $\RR^2$. This particular isomorphism was chosen so that simple closed curves, thought of as measured foliations, map to integer points. The numbers $(a,b,c)$ are called the {\em triangle lengths} and $(x,y)=\pi(a,b,c)$ are the {\em collar parameters}. It is routine to check that the map $\pi_{\Hcal}=\pi|{\Hcal}:\Hcal\rightarrow\RR^2$ is a homeomorphism, so the collar parameters determine the triangle lengths.

A collar parameter $p\in\RR^2$ gives rise to both a measured foliation and a hyperbolic metric on $\stdan$, and thus also on the universal cover $\widetilde \stdan$. 

\vspace{3mm}

%%%%%%%%%%%%%%%%%%%%%%%%%%%%%%%%%%%%%%%%%%%%%%%%%%%%%%%%%%%%%%%

 \subsection{Parameterizing Teichm\"uller space}

%%%%%%%%%%%%%%%%%%%%%%%%%%%%%%%%%%%%%%%%%%%%%%%%%%%%%%%%%%%%%%%

In the following $\Sigma$ is a closed orientable connected surface of genus $g\ge 2$ and, as in \S \ref{ss:ShortsDecomp},
$$\Acal=\{\alpha_i:1\le i\le 3g-3\}$$ 
is a set of disjoint simple closed curves in $\Sigma$ such that the closure of each component of $\Sigma\setminus\Acal$ is a topological pair of pants. 
A point in $\Teich(\Sigma)$ is uniquely determined by collar parameters $(x_i, y_i) \in \RR^2$ for each $\alpha_i$ as follows.

Choose  a set of pairwise disjoint annuli 
$\{ A_i :1\le i\le 3g-3\}$
such that  $A_i$ contains $\alpha_i$ in the interior for each $i$. Choose an identification of $A_i$ and the standard annulus $\stdan$. The collar parameters $(x_i,y_i)$ determine a hyperbolic metric on $A_i$, using the triangle lengths $\pi^{-1}(x_i,y_i)$.
 
The closure of a component of $\Sigma\setminus\Acal$ is a hyperbolic pair of pants. The hyperbolic structure on each pair of pants is determined by the collar parameters. The basepoints of the standard collars are on the seams of the pants they are in. This determines how standard collars are glued to pants and thus determines a hyperbolic metric on $\Sigma$ and a parameterization
$$\eta: \prod\RR^2\longrightarrow\Teich(\Sigma).$$

\vspace{3mm}

%%%%%%%%%%%%%%%%%%%%%%%%%%%%%%%%%%%%%%%%%%%%%%%%%%%%%%%%%%%%%%%
 
\subsection{Parameterizing measured foliations}\label{ss:Foliations}

%%%%%%%%%%%%%%%%%%%%%%%%%%%%%%%%%%%%%%%%%%%%%%%%%%%%%%%%%%%%%%%

Our point of view is that a measured foliation $|\omega|$ is determined by a $1$-form $\omega$. 
We regard two measured foliations as {\em equivalent}
if they determine the same length functions on $\Scal$.
If $\omega=0$ on a subsurface then
the foliation on that subsurface 
is not important. 
The discussion below follows the terminology of \cite[Section 6.2]{FLP}.  We wish to concentrate the transverse measure in the annuli.
We will define a {\em standard measured foliation}
on shorts and on an annulus. Fix a shorts decomposition of the surface $\Sigma$. Then a measured foliation on $\Sigma$ is {\em standard}
if the restriction to each pair of shorts and each  annulus is standard. It follows from \cite{FLP} that:
\begin{proposition}\label{standardfolation} Every measured foliation on a surface is  equivalent to a standard one.\end{proposition}

\begin{center}
\begin{figure}[H]
\includegraphics[height=1.15in]{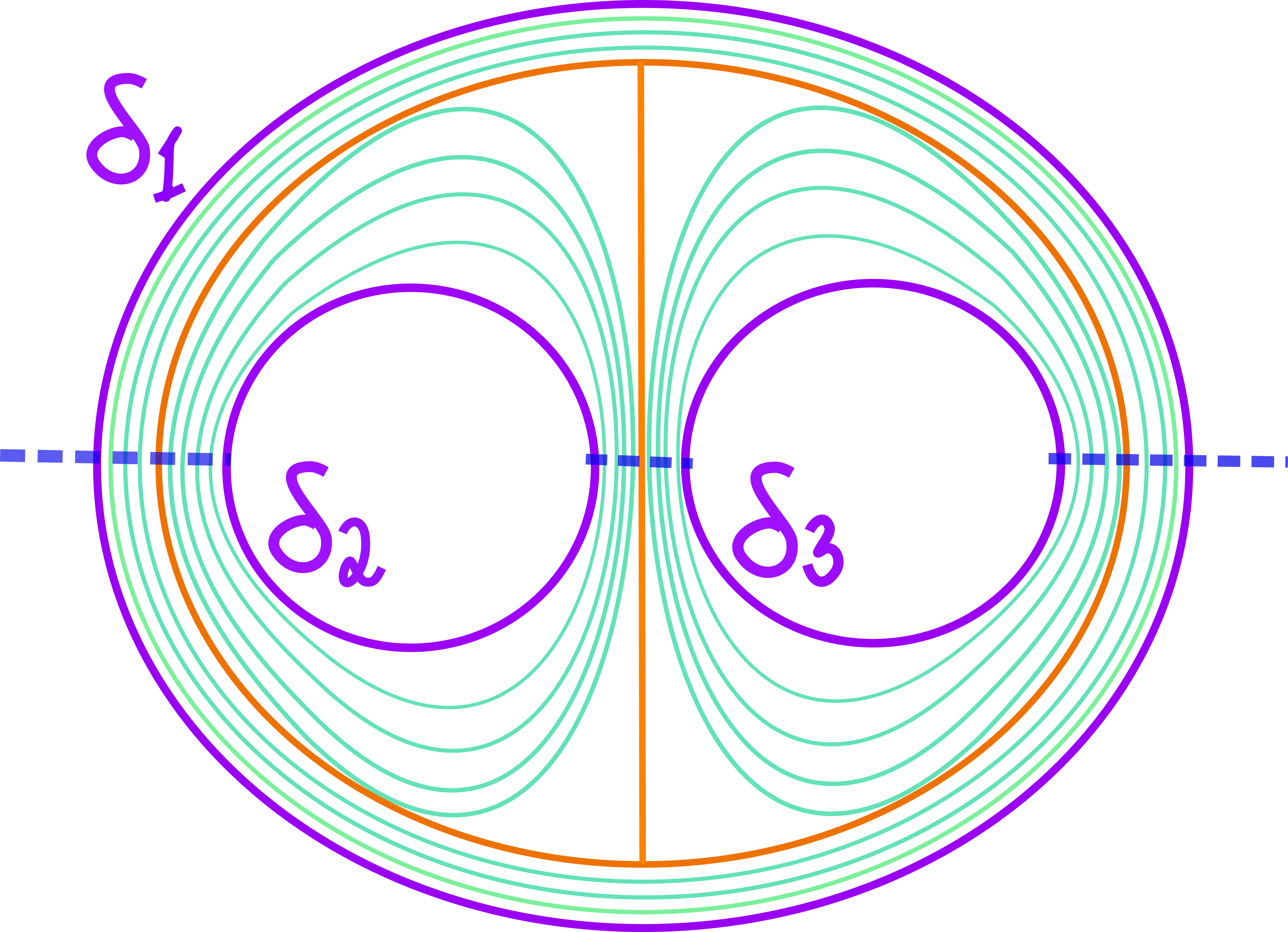}\quad\quad\quad
\includegraphics[height=1.15in]{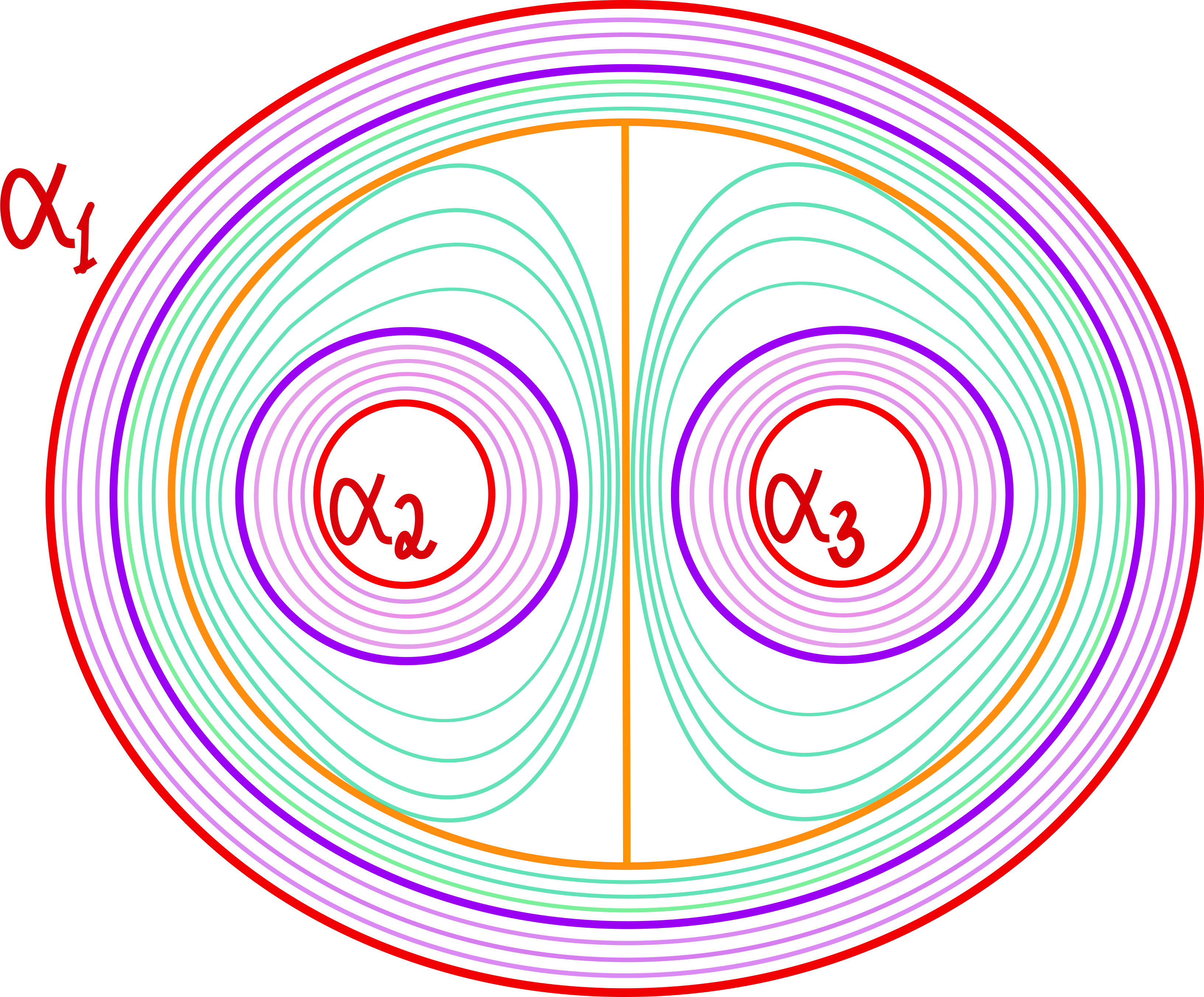}
\vspace*{-2mm}
\caption{The measured foliation in the case of $m_1=m_2=m_3=0$ on shorts (left) and pants (rights). The line of reflection in the shorts is the blue dotted line.} \label{zeropants}
\end{figure}
\end{center}

\begin{center}
\begin{figure}[H]
\includegraphics[width=6in]{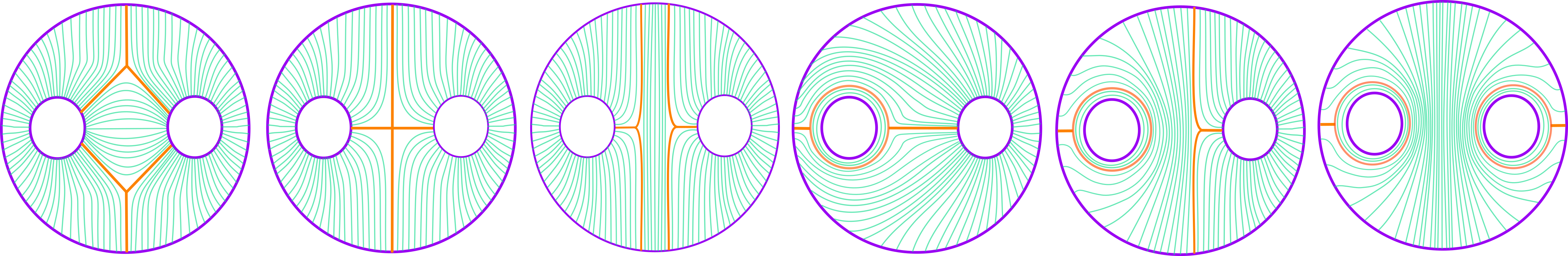}
\caption{Possible foliations on the shorts \cite[Figure 6.6]{FLP}: from left to right these are where $m_1+m_2+m_3 > 2max\{m_1,m_2,m_3\}$,  and
$m_1 = m_2+m_3$,  and
$m_1 > m_2 + m_3$, and
$m_1=m_3, m_2=0$, and
$m_1>m_3, m_2=0$, and
$m_1>0, M_2=m_3=0$.} \label{foliatedshorts}
\end{figure}
\end{center}

Let $P$ be a pair of shorts with  boundary components $\delta_1,\delta_2,\delta_3$ (in \cite{FLP} the corresponding boundary components are called $\gamma_1,\gamma_2,\gamma_3$). Given $m_1,m_2,m_2\ge 0$ we
define  a {\em standard measured foliation} on  $P$ such that
the transverse measure of $\delta_i$ is $m_i$. Except for the case of $(m_1,m_2,m_3)=(0,0,0)$ the leaves of the foliation are shown in \cite[Figure 6.6]{FLP},
 but modified as follows. If $m_i=0$ then an annulus neigborhood of  $\delta_i$ is foliated by smooth circles, and with transverse measure zero. The result
 is shown in  Figure \ref{foliatedshorts}.
The remaining case of $(m_1,m_2,m_3)=(0,0,0)$ is shown in Figure \ref{zeropants}.
\vspace*{-9mm}

 In each case there is an automorphism of $P$ that is a reflection that fixes  the union of three arcs, one connecting
each pair of boundary components, and
that is measure preserving.

\vspace*{-8mm}

\noindent  

It remains to define a standard measured foliation on an annulus. They are given by a linear $1$-form. There are  two kinds, depending on whether the leaves are circles,
or arcs connecting the two boundary components. In the first case the transverse measure might be zero. 

\begin{definition}[triangle equality]
A point $x\in\RR^3$ {\em satisfies the  triangle equality} if it is in the set
\begin{equation}\label{e:TriangleEquality}
\Ccal=\{(a,b,c)\ :\ a+b+c=2\max\{a,b,c\}\quad \&\quad a,b,c\ge 0\}.
\end{equation}
\end{definition}

This is the cone from $0$ on a $2$-simplex. Figure \ref{MagicPic} shows how the subset $\Hcal$ of $\RR^3$ sits inside $\Ccal$ like a hyperboloid inside its lightcone: they are asymptotic at infinity. Moreover $\pi|:\Delta\rightarrow\RR^2$ is a homeomorphism.

A measured foliation $\mu$ on $\stdan$ is  {\em linear} if it is covered by a measured foliation $|df|$ on $\widetilde \stdan\subset\RR^2$ given by the restriction of some linear map $f:\RR^2\to \RR$. A Euclidean line segment in $\stdan$ is either transverse to the foliation, or else contained in a leaf. We assign {\em lengths} $(a,b,c)$ to the sides of the reference triangle in $\stdan$ by integrating $|df|$ along each side. Then $(a,b,c)\in\Ccal$ and $h=a$. We again refer to $\pi(a,b,c)\in\RR^2$ as {\em CP parameters}, and
 they determine these lengths for measured foliations.
  
  \parpic[r]{\includegraphics[width=1in]{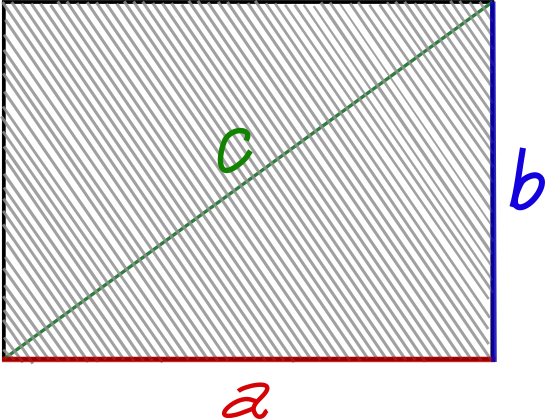}}
Given edge lengths $v=(a,b,c)\in\Ccal$, if $c\ge \max(a,b)$ then $c=a+b$, otherwise $c=|a-b|$. Let $\mf(v)$ be the linear measured foliation on $\stdan$ that assigns lengths $a,b,c$ to the standard unit
vectors $e_1,e_2,e_1+e_2$ respectively. Then $\mf(v)$ lifts to a measured foliation $|\omega_v|$ where $\omega_v$ is the $1$-form on $\RR^2$ given by
\begin{equation}
\omega_v=
\left\{
\begin{array}{lll}
a.dx+b.dy &  & c=a+b\\
a.dx-b.dy & & c=|a-b|
\end{array}\right.
\end{equation}
The figure above shows a foliation on $\RR^2$ such that $\omega_v$ vanishes on the tangent spaces of the leaves.
When $v=0$ we define the leaves to be the circles given by $\ker dy$.

An assignment of a point in $\pi(\Delta)=\RR^2$ to each annulus in an annulus-shorts decomposition of $\Sigma$  determines a measured foliation on $\Sigma$ and a parameterization
$$\mu:\prod\RR^2\longrightarrow\MF(\Sigma).$$

\vspace{3mm}

%%%%%%%%%%%%%%%%%%%%%%%%%%%%%%%%%%%%%%%%%%%%%%%%%%%%%%%%%%%%%%%

\subsection{CP maps}

%%%%%%%%%%%%%%%%%%%%%%%%%%%%%%%%%%%%%%%%%%%%%%%%%%%%%%%%%%%%%%%

%\begin{definition}[collar parameter (cp) coordinates, foliation map, hyperbolization map] 
The {\em collar parameters} on Teichm\"uller space and measured foliations are the maps $\Theta_{\Teich}= \eta^{-1}: \Teich(\Sigma) \to\RR^{6g-6}$ and $\Theta_{\MF}=\mu^{-1}:\MF(\Sigma)\to\RR^{6g-6}$. The map $m=\mu\circ \eta^{-1}:\Teich(\Sigma)\to\MF(\Sigma)$ is called the {\em  foliation map}, and the inverse $h: \MF(\Sigma) \to\Teich(\Sigma)$ the {\em hyperbolization map}.
%\end{definition}

\vspace{5mm}

%%%%%%%%%%%%%%%%%%%%%%%%%%%%%%%%%%%%%%%%%%%%%%%%%%%%%%%%%%%%%%%
 
\section{The Thurston compactification is the Radial Compactification}

%%%%%%%%%%%%%%%%%%%%%%%%%%%%%%%%%%%%%%%%%%%%%%%%%%%%%%%%%%%%%%%

Suppose that $ds$ is a positive semi-definite quadratic form on $\Sigma$. There is a {\em length function}
$$\Lcal(ds):\Scal\rightarrow\RR$$ defined as follows. 
Given an element $\sigma$ of $\Scal$ then $$(\Lcal(ds))(\sigma)=\inf\int_\gamma ds$$ where the infimum is taken over all
simple closed curves $\gamma$ in the isotopy class $\sigma$. We are interested in applying the length function to $ds$, when it is given either by a hyperbolic metric, or by a transversally measured foliation on $\Sigma$.

Choose a finite symmetric generating set $W\subset \pi_1(\Sigma)$. Set $W^1=W$ and $W^{n+1}=\{x.y:x\in W, y\in W^n\}$. For $g\in\pi_1\Sigma$ define $w:\pi_1\Sigma\to Z$ by $w(g)=\min\{n:\exists h\ \  hgh^{-1}\in W^n\}$. The number $w(g)$ is called the {\em conjugacy word length} of $g$ and is the minimum length of a word in the elements of $W$ that is conjugate to $g$. The hyperbolic structure, $\eta(\sigma)$, and the measured foliation, $\mu(\sigma)$, give length functions on $\Scal$ that differ by less than a fixed multiple of word length:

\begin{mainthmA}\label{thm1} Given a conjugacy word length $w$ on $\pi_1\Sigma$, there is a constant $C>0$ such that $m:\Teich(\Sigma)\to\MF(\Sigma)$ satisfies:
$$ \forall\ \rho\in\Teich(\Sigma)\qquad |\ \Lcal(\rho)-\Lcal(m(\rho))\ |\le C\cdot w$$
\end{mainthmA}

Using the embedding $\RR^n\hookrightarrow\RR^n$ given by $v\mapsto v/(1+\|v\|)$ the {\em radial compactification} of $\RR^n$ is the unit ball $B=\{v\in\RR^n:\|v\|\le 1\}$. The {\em Thurston compactification} is $\overline\Teich(\Sigma)=\Teich(\Sigma)\sqcup{\mathbb P}(\MF)$.
Since $\mu(t v)=t\cdot \mu(v)$ for $t>0$, it follows from Theorem A that the length functions $t^{-1}\Lcal(\eta(t v))$ converge to the length function of $\mu(v)$ provided $v\ne 0$.

\begin{mainthmB} Using $\Theta_{\Teich}$ coordinates to identify $\Teich(\Sigma)\equiv\RR^{6g-6}$, the Thurston compactification is the radial compactification of $\RR^{6g-6}$. For $0\ne v\in\RR^{6g-6}$, we have $\lim_{t\to\infty}\eta(t v)=[\mu(v)]\in \overline\Teich(\Sigma)$.
\end{mainthmB}

\vspace{3mm}

%%%%%%%%%%%%%%%%%%%%%%%%%%%%%%%%%%%%%%%%%%%%%%%%%%%%%%%%%%%%%%%

\section{Realizing the compactification with quadratic differentials.}

%%%%%%%%%%%%%%%%%%%%%%%%%%%%%%%%%%%%%%%%%%%%%%%%%%%%%%%%%%%%%%%

Suppose $ds_0$ is some Riemannian metric on $\Sigma$, not necessarily hyperbolic, called the {\em background metric}. We show that, {\em after a suitable isotopy}, a hyperbolic metric  on $\Sigma$ differs from some measured foliation on $\Sigma$ by less than a fixed multiple of $ds_0$. Then integration along geodesics shows the two length functions are close, provided these geodesics are not too long in the background metric. We formalize this  with the following.

\begin{definition} A  seminorm $ds$ on $\Sigma$ is {\em $C$-efficient with respect to the background metric $ds_0$} if  for each $g\in\pi_1\Sigma$ there is a $ds$-geodesic $\alpha:S^1\to\Sigma$ that is freely homotopic to a loop representing $g$, and $\ell(\alpha,ds_0)\le C\cdot w(g)$. A set of seminorms $\Ncal$ is {\em uniformly efficient} if there is $C>0$ such that all the seminorms in $\Ncal$ are $C$-efficient.
\end{definition}

Since $\Sigma$ is compact any two background metrics are bilipschitz. Thus whether or not a set of seminorms is uniformly efficient does not depend on the choice of background metric.
  
The space of quadratic differentials $\Qcal(\Sigma)$ on $\Sigma$ contains the subspace, $\widetilde{\Teich}(\Sigma)$, of hyperbolic metrics and the subspace, $\widetilde{\MF}(\Sigma)$, of measured foliations. There are natural projections
$\pi_{\Teich}:\widetilde{\Teich}(\Sigma)\to\Teich(\Sigma)$ and $\pi_{\MF}:\widetilde{\MF}(\Sigma)\to\MF(\Sigma)$, and $\pi_{\Teich}$ is a fiber bundle with fiber the group of diffeomorphisms isotopic to the identity. These maps have sections:

\begin{mainthmC}[Efficient Realization Theorem]\label{thm2} 
Suppose  $\Sigma$ is a closed orientable surface with genus at least $2$. Then there are embeddings
$\widetilde{\eta}:\RR^{6g-6}\to\widetilde\Teich(\Sigma)$ and 
$\widetilde{\mu}:\RR^{6g-6}\to\widetilde{\MF}(\Sigma)$ 
such that $\mu=\pi_{\MF}\circ\widetilde{\mu}$ and $\eta=\pi_{\Teich}\circ\widetilde{\eta}$, with $\widetilde{\eta}(t x)=t\widetilde{\eta}(x)$ for each $t\ge0$. Moreover, given a background  metric $ds_0$ on $\Sigma$,
there is $C=C(ds_0)>0$  so that the image of $\widetilde{\mu}$ and of $\widetilde{\eta}$ are uniformly efficient and 
\begin{equation}\label{maineq}\forall x\in\RR^{6g-6}\qquad
  |\ \widetilde{\eta}(x)-
 \widetilde{\mu}(x)\ |\le C\cdot  |ds_0| 
 \end{equation}
\end{mainthmC}
To define $\widetilde\mu$ and $\widetilde
\eta$ involves writing down explicit metrics and measured foliations on shorts that match standard metrics on collars along the boundary. 
 
\vspace{3mm}

%%%%%%%%%%%%%%%%%%%%%%%%%%%%%%%%%%%%%%%%%%%%%%%%%%%%%%%%%%%%%%%
 
\section{Converting between Collar Parameters and Fenchel-Nielsen coordinates}
 
%%%%%%%%%%%%%%%%%%%%%%%%%%%%%%%%%%%%%%%%%%%%%%%%%%%%%%%%%%%%%%%

We provide here the coordinate change maps between the CP coordinates we have defined, and the classical Fenchel-Nielsen coordinates on Teichm\"uller space
and Dehn-Thurston coordinates for measured foliations. 

\begin{proposition}\label{FNDT}
Suppose $(2\ell,2\tau)\in\RR^+\times\RR$ are Fenchel-Nielsen coordinates \cite{Wolpert}. The triangle lengths  $(a,b,c)$ are given by
\begin{equation}
\begin{array}{rcl}
a & = &\ell\\ 
b & = &\cosh^{-1}(\cosh \tau\coth \ell)\\  
c & = &\cosh^{-1}(\cosh(\ell-\tau)\coth \ell)
\end{array}
\end{equation}
and the collar parameters are given by
\begin{equation}
\begin{array}{rcl}
x & = &4\ell-2\cosh^{-1}\left(\cosh \tau\coth \ell \right)-2\cosh^{-1}\left(\cosh(\ell-\tau)\coth \ell \right)\\
y & = &2\cosh^{-1}(\cosh \tau\coth \ell )-2\cosh^{-1}(\cosh(\ell-\tau)\coth \ell )
\end{array}
\end{equation}
If $(2\ell,2\tau)\in\RR_{\ge 0}\times\RR$ are Dehn-Thurston coordinates \cite{Feng} then  the triangle lengths $(a,b,c)$  are given by $a=\ell$  and $b=|\tau|$ and $c=|\ell-\tau|$ and the collar parameters by
\begin{equation}
\begin{array}{rcl}
x & = & \ 4\ell-2|\tau|-2|\ell-\tau|\\
y & = & \ 2|\tau|-2|\ell-\tau|
\end{array}
\end{equation}
\end{proposition}
When $\ell$ is large then $\coth\ell\approx 1.$ Observe that replacing $\coth\ell$ by $1$ in (6) yields (7).
\vspace{3mm}

%%%%%%%%%%%%%%%%%%%%%%%%%%%%%%%%%%%%%%%%%%%%%%%%%%%%%%%%%%%%%%%
 
\section{The collar equation and Teichm\"uller space of a once-punctured torus}

%%%%%%%%%%%%%%%%%%%%%%%%%%%%%%%%%%%%%%%%%%%%%%%%%%%%%%%%%%%%%%%

The collar equation is also the equation of the character variety for the Teichm\"uller space of finite area hyperbolic structures on a once-punctured torus $T$. This follows because the worst case for a standard collar is given by $T$, where
  there is a single self intersection point on the boundary of a standard collar. This point determines a reference triangle in $T$ and this triangle determines the metric on $T$
  up to isotopy.  Here are the details.

\begin{center}
\begin{figure}[h]
\includegraphics[scale=.27]{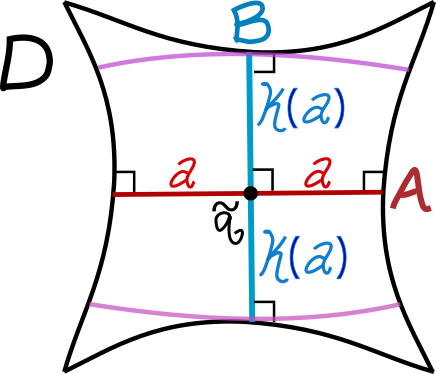}
        \hspace{15mm}
\includegraphics[scale=.29]{LabeledTorusNew.png}
 \caption{ (a) The left-hand image depicts an 
ideal quadrilateral $D$ in $\HH^2$ such that the common perpendiculars, $A$ and $B$, to opposite sides of $D$ are orthogonal. If the length of $A$ is $2a$, then the length of $B$ is $2\kappa(a)$. The $2\kappa(a)$-hypercycles are depicted in purple. (b) The right-hand image depicts the punctured torus obtained by identifying the opposite sides of $D$ using isometries that translate along the common perpendiculars. $\alpha$ is the image of $A$ and $\beta$ is the image of $B$. The image of the $2\kappa(a)$-hypercycles correspond to the standard annulus boundary image.} \label{f:TorusTeich}
\end{figure}
\end{center}

\begin{center}
\begin{figure}[H]
\includegraphics[scale=.4]{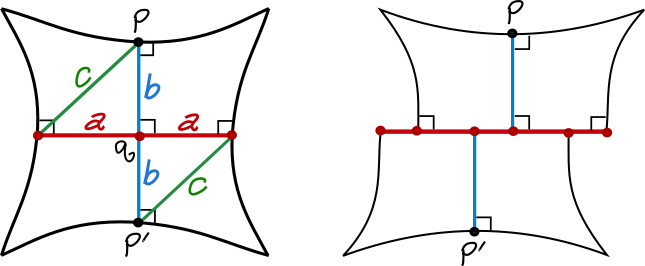}
\vspace*{3mm}
 \caption{$\tr[A,B]=-2\quad\Rightarrow\quad\cosh^2a+\cosh^2b+\cosh^2c=2\cosh a\cosh b\cosh c$} \label{Fig10}
\end{figure}
\end{center}

Refer to Figures \ref{f:TorusTeich} and \ref{Fig10}. Let $D$ be an ideal quadrilateral in $\HH^2$ such that the common perpendiculars, $A$ and $B$, to opposite sides of $D$ are orthogonal. Then there is an $a\in \R_{>0}$ so that the lengths of these common perpendiculars are $2a$ and $2\kappa(a)$. Now $D$ is a fundamental domain for a hyperbolic metric on $T$. This is obtained by identifying the opposite sides of $D$ using isometries that translate along the common perpendiculars.
The image of $A$ in $T$ is a simple closed geodesic $\alpha$. The standard collar of $\alpha$ meets itself at one point $p$, that is the image of the endpoints of $A$. The image of $B$ is a closed geodesic $\beta$ on $T$ that is orthogonal to $\alpha$ at the point $q$. There is a third closed geodesic $\gamma$ on $T$ that contains $p$ and $q$ and is homotopic to $\alpha\cdot \beta$. Then the reference triangle has side lengths $a,b,c$ that are the half-lengths of the geodesics $\alpha$, $\beta$, and $\gamma$.

\vspace{-5mm}

If $A,B\in\SL(2,\RR)$ are the holonomies of $\alpha$ and $\beta$ then $C=AB$ is the holonomy of $\gamma$. It follows from the trace relation that
\begin{equation}\label{trreln}
\tr[A,B]=-2=(\tr A)^2+(\tr B)^2+(\tr C)^2-\tr A\tr B\tr C-2.
\end{equation}
The relationship $\tr A=2\cosh a$ between the trace and translation length yields the collar equation in this case.

The hyperbolic structure on $T$ has the property that $\alpha$ and $\beta$ are orthogonal. Any finite area structure on $T$ can be obtained from some such $\alpha$ and $\beta$ by an earthquake \`a la Fenchel-Nielsen along $\alpha$ by some distance. One can picture
the triangulation on the resulting structure by cutting $D$ along $A$ and sliding the bottom half sideways, see Figure \ref{Fig10}. It follows that the standard collar of $\alpha$ in the resulting structure still has a single point of self intersection.
The geodesic $B$ is replaced by the geodesic connecting $p$ and the image $p'$ of $p$ under the sideways slide.

\vspace{-10mm}

\begin{proposition} 
Let $T$ be a once-punctured torus and $\alpha$, $\beta$ a generating set for $\pi_1(T)$. Let $\Tcal(T)$ denote the Teichm\"uller space of finite area hyperbolic metrics on $T$. Then there exists a homeomorphism $\theta:\Tcal(T)\longrightarrow\Hcal$ so that $\theta(\rho)=(a,b,c)$ 
are half the lengths of geodesic representatives of $\alpha$, $\beta$, and $\alpha.\beta$ respectively.
\end{proposition}

\vspace{5mm}
\vspace{5mm}

\small

\bibliography{Coords} 

\bibliographystyle{amsplain}

\end{document}